\documentclass[12pt]{article}
\usepackage{graphicx,psfrag,epsfig}
\usepackage{amssymb,amsmath,amscd,amsthm}
\usepackage{graphicx,psfrag,epsfig}

\usepackage{graphicx}
\usepackage[active]{srcltx}

\newtheorem{theorem}{Theorem}[section]

\newtheorem{lemma}[theorem]{Lemma}

\date{}

\begin{document}

\date{}
\title{Branching
diffusion in inhomogeneous media}
\author{L.
Koralov\footnote{Dept of Mathematics, University of Maryland,
College Park, MD 20742, koralov@math.umd.edu}
} \maketitle
\begin{abstract}
We investigate the long-time evolution of branching diffusion
processes (starting with a finite number of particles) in
inhomogeneous media. The qualitative behavior of the processes
depends on the intensity of the branching. In the super-critical
regime, we describe the asymptotics of the number of particles in
a given domain. In the sub-critical and critical regimes, we show
that the limiting number of particles is finite and describe its
distribution.
\end{abstract}

\section{Introduction}

 Consider a collection
of particles in $ \mathbb{R}^d$ that move diffusively and
independently. Besides the diffusive motion, the particles can
duplicate with the rate of duplication $\beta v(x)$, $x \in
\mathbb{R}^d$, where $x$ is the position of a given particle, $v$
is a continuous non-negative compactly supported function and
$\beta \geq 0$ is a parameter controlling the duplication rate.
Both copies start moving independently immediately after the
duplication.

There is a number of results describing the processes in the case
when $v$ in concentrated in one point. See, for example,
\cite{DaFl}, \cite{DFLM}, \cite{DFL}, \cite{FL} for the study of
superprocesses and \cite{ABY}, \cite{BY1}, \cite{ABY2}, \cite{Ya}
for the asymptotic properties of branching random walks. We now
will give a detailed description of the behavior of the branching
diffusion in $ \mathbb{R}^d$ in the case when $v$ is a compactly
supported function. Our main goal is to describe the distribution
of particles when $t$ is large. The asymptotics depends on whether
$\beta$ is above, at, or below the critical value $\beta_{\rm
cr}$, which is the infimum of values of $\beta$ for which the
operator
\begin{equation} \label{opera}
 L^\beta u(x) = \frac{1}{2} \Delta u(x) +
 \beta v(x) u(x)
\end{equation}
  has a
positive eigenvalue. This is the operator in the right hand side
of the equations on the particle density and higher order
correlation functions, given below.

 We will show that  for $\beta >
\beta_{\rm cr}$ the number of particles in a given region $U$ at
time $t$ has the asymptotics
\[
n_t(U) \sim e^{\lambda_0(\beta)t }\xi \int_U \varphi(y) d y,
\]
where $\lambda_0(\beta)$ is the largest eigenvalue of $L^\beta$,
$\xi$ is a random variable that depends on the initial
configuration of particles (assumed to be finite in number) and
$\varphi$ is a deterministic function (limiting density profile).
Intuitively, the presence of the random variable $\xi$ reflects
the effect of branching at random times while the number of
particles is small. After the number of particles becomes
sufficiently large, it keeps growing nearly deterministically due
to the fact that the bulk of particles is located near the support
of $v$ and due an `averaging' effect in the branching mechanism.

When $\beta < \beta_{\rm cr}$  (and $d \geq 3$), the effect of
transience of the diffusion outweighs the branching, and all the
particles will eventually wander off to infinity (assuming that
initially there were finitely many particles). It will be shown
that for $\beta < \beta_{\rm cr}$ the total number of particles
tends, as $t \rightarrow \infty$, to a finite random limit, whose
distribution will be identified.

The case when $\beta = \beta_{\rm cr}$ is interesting when $d \geq
3$ ($\beta_{\rm cr} = 0$ for $d =1,2$). In this case the total
number of particles will be shown to tend to a finite limit almost
surely, although the expectation of the total number of particles
tends to infinity.

Some of the corresponding results for the processes on the lattice
with branching at the origin were obtained in \cite{ABY},
\cite{BY1}, \cite{ABY2}, \cite{Ya} . The main difficulty here,
compared to the latter series of papers, is that the explicit
formulas for the resolvent of the generator that were helpful in
analyzing the processes with branching at the origin are not
available now. Besides allowing the treatment of the general
potential $v$, the techniques developed in this paper will allow
us to study some of the more intricate properties of the limiting
distribution: the fluctuations of the local number of particles
conditioned on the total number of particles in the super-critical
case, the growth of the region containing the particles and the
distribution of the number of particles in the regime of large
deviations (near the edge of the region containing the particles).
These properties will be the subject of a subsequent paper.

The results on the large time asymptotics (Sections
\ref{scc}-\ref{critc}) will be obtained from the equations on the
particle density and higher order correlation functions derived in
Section~\ref{momeq}. The analysis will be based on the spectral
representation of solutions in the appropriate function spaces
followed by the asymptotic analysis of integrals with integrands
that depend on several parameters. Some of the techniques are
related to those employed by us in the study of a polymer
distribution in a mean field model (\cite{CKMV1}). There, the
asymptotics of a single equation (rather than a recursive system
of equations) in the parameters $t$ and $\beta$ was examined.

Finally, let us mention a number of recent papers on the parabolic
Anderson model, where $v$ is a stationary random field (see, for
example, \cite{GM1}, \cite{GM2}, \cite{CaMo}, \cite{GKM},
\cite{GK}, \cite{KLMS}). When $v$ is random, the behavior of the
solution  to (\ref{fmo})-(\ref{ico}) essentially depends on nature
of the tails of the distribution of $v$. It has been shown that in
many cases the main contribution to the creation of the total
number of particles is given by the isolated high peaks of the
random potential. Moreover, when $t$ is large, the bulk of the
solution is located near one of those peaks with high probability.
This adds to the importance of the study of branching diffusions
in the case when $v$ is localized.

\section{Equations on correlation functions}
\label{momeq}

 Let
$B_\delta$ be a ball of radius $\delta$ in $ \mathbb{R}^d$. For $t
> 0$ and $x,y_1,y_2,... \in \mathbb{R}^d$ with all $y_i$  distinct,
define the particle density $\rho_1(t,x,y_1)$ and the higher order
correlation functions $\rho_n(t,x,y_1,...,y_n)$ as the limits of
probabilities of finding $n$ distinct particles in
$B_\delta(y_1)$,...,$B_\delta(y_n)$, respectively, divided by
${\rm Vol}^n (B_\delta)$, under the condition that there is a
unique particle at $t = 0$ located at $x$. We extend
$\rho_n(t,x,y_1,...,y_n)$ by continuity to allow for $y_i$ which
are not necessarily distinct. For fixed $y_1$, the density
satisfies the equation
\begin{equation} \label{fmo}
\partial_t \rho_1(t,x,y_1) = \frac{1}{2} \Delta \rho_1(t,x,y_1) +
\beta v(x) \rho_1(t,x,y_1),
\end{equation}
\begin{equation} \label{ico}
\rho_1(0,x,y_1) = \delta_{y_1}(x).
\end{equation}
Indeed, let $s,t > 0$. Then we can write
\begin{equation} \label{appro}
\rho_1(s+t,x,y_1) = (2 \pi s)^{-\frac{d}{2}} \int_{ \mathbb{R}^d}
e^{\frac{|x - z|^2}{2s}} \rho_1(t,z,y_1)d z + \beta v(x) s
\rho_1(t,x,y_1) + \alpha(s,t,x,y_1),
\end{equation}
 where the term with the
integral on the right hand side is due to the effect of the
diffusion on the interval $[0,s]$, the second term is due to the
probability of branching on $[0,s]$, and $\alpha$ is the
correction term. The correction term is present since (a) more
than one instance of branching may occur before time $s$, and (b)
even if a single branching occurs between the times $0$ and $s$,
then the original particle will be located not at $x$ but at a
nearby point and the intensity of branching there is slightly
different from $\beta v(x)$. It is clear that $\lim_{s \downarrow
0} \sup_{x,y \in \mathbb{R}^d} \alpha(s,t,x,y)/s = 0$. After
subtracting $\rho_1(t,x,y_1)$ from both sides of~(\ref{appro}),
dividing by $s$ and taking the limit as $s \downarrow 0$, we
obtain (\ref{fmo}).

The equations on $\rho_n$, $n > 1$, are somewhat more complicated:
\begin{equation} \label{manyp}
\partial_t \rho_n(t,x,y_1,...,y_n) = \frac{1}{2} \Delta \rho_n(t,x,y_1,...,y_n) +
\beta v(x) \left( \rho_n(t,x,y_1,...,y_n) +  H_n (t,x,y_1,...,y_n)
\right) ,
\end{equation}
\begin{equation} \label{initu}
\rho_n(0,x,y_1,...,y_n) \equiv 0.
\end{equation}
Here
\[
H_n (t,x,y_1,...,y_n) = \sum_{U \subset Y, U \neq \emptyset}
\rho_{|U|} (t,x,U) \rho_{n-|U|}(t,x, Y\setminus U),
\]
where $Y = (y_1,...,y_n)$, $U$ is a proper non-empty subsequence
of $Y$, and $|U|$ is the number of elements in this subsequence.
Equation (\ref{manyp}) is derived similarly to (\ref{fmo}). The
combinatorial term $H_n$ appears after taking into account the
event that there is a single branching on the time interval
$[0,s]$, the descendants of the first particle are found at the
points in $U$ at time $s+t$, while the descendants of the second
particle are found at the points of $ Y \setminus U$, with the
summation over all possible choices of $U$.

\section{Analytic Properties of the Resolvent} \label{sectt}

Here we recall some basic facts about the operator
$L^\beta:L^{2}({\mathbb{R}}^{d})\rightarrow
L^{2}({\mathbb{R}}^{d})$ (see (\ref{opera}))  and its resolvent
$R^{\beta }_\lambda =(L^{\beta }-\lambda )^{-1}$. We will assume
that $v \geq 0$ is continuous, compactly supported and not
identically equal to zero. It is
well-known that the spectrum of $L^{\beta }$ consists of the absolutely continuous part $%
(-\infty ,0]$ and at most a finite number of non-negative
eigenvalues:
\[
\sigma (L^{\beta })=(-\infty ,0]\cup \{\lambda _{j}\},\text{ \ \
}0\leq j\leq N,\text{ \ \ }\lambda _{j}=\lambda _{j}(\beta )\geq
0.
\]
We enumerate the eigenvalues in the decreasing order. Thus, if
$\{\lambda _{j}\}\neq \emptyset $, then $\lambda _{0}=\max \lambda
_{j}$.  Thus the resolvent $R^{\beta }_\lambda:
L^{2}({\mathbb{R}}^{d})\rightarrow L^{2}({\mathbb{R}}^{d})$ is a
meromorphic operator valued function on $\mathbb{C}^{\prime }
 =\mathbb{C}\backslash (-\infty ,0]$.
\begin{lemma}
\label{lmonot} There exists $\beta _{cr} \geq 0$ (which will be
called the critical value of $\beta $) such that $\sup \sigma
(L^{\beta })=0$ for $\beta \leq \beta _{cr}$ and $\sup \sigma
(L^{\beta })=\lambda _{0}(\beta)>0$ for $\beta >\beta _{cr}$. For
$\beta
>\beta _{cr}$ the eigenvalue $\lambda _{0}(\beta )$ is a strictly
increasing and continuous function of $\beta $. Moreover,
$\lim_{\beta \downarrow \beta _{cr}}\lambda (\beta )=0$ and
$\lim_{\beta \uparrow \infty
}\lambda (\beta )=\infty $. 
\end{lemma}

The proof of this lemma is standard (see Lemma 4.1 of
\cite{CKMV1}).  Denote the kernel of $R^{\beta }_\lambda$ by
$R^{\beta }_\lambda(x,y).$ If $\beta =0$, the kernel depends on
the difference $x-y$ and will intermittently use the notations
$R^{0}_\lambda(x,y)$ and $R^{0}_\lambda(x-y)$. The kernel $R^{0
}_\lambda(x)$ can be expressed through the Hankel function $%
H_{\nu }^{(1)}$:
\begin{equation}  \label{ha2}
R^{0}_\lambda(x)=
c_dk^{d-2}(k|x|)^{1-\frac{d}{2}}H_{\frac{d}{2}-1}^{(1)}(i\sqrt2k|x|),~~
k=\sqrt{\lambda},~~\mathrm{Re}k > 0.
\end{equation}

 We shall say that  $f \in C_{\rm \exp}( \mathbb{R}^d)$ (or simply
 $C_{\rm \exp}$)
if $f$ is continuous and
\[
||f||_{C_{\rm \exp}( \mathbb{R}^d)} = \sup_{x \in \mathbb{R}^d}(
|f(x)| e^{{|x|^2}}) < \infty.
\]
The space of bounded continuous functions on $\mathbb{R}^d$ will
be denoted  by $C( \mathbb{R}^d)$ or simply $C$. The following
lemma will be proved in the Appendix.
\begin{lemma} \label{rbetaW}
The operator $R^\beta_\lambda: C_{\exp}( \mathbb{R}^d) \rightarrow
C( \mathbb{R}^d)$ is meromorphic in $\lambda \in \mathbb{C}'$. Its
poles are of the first order and are located at eigenvalues of the
operator $L^{\beta}$. For each $\varepsilon
> 0$ and some $\Lambda=\Lambda(\beta)$, the operator is uniformly bounded
in $\lambda \in \mathbb{C}'$,
 $|{\rm arg} \lambda| \leq \pi -
\varepsilon$, $|\lambda| \geq \Lambda$. It is of order
$O(1/|\lambda|)$ as $\lambda \rightarrow \infty$, $|{\rm arg}
\lambda| \leq \pi - \varepsilon$.

If $d \geq 3$ and $\beta \in [0, \beta_{\rm cr})$, then
$R^\beta_\lambda$ has the following asymptotic behavior as
$\lambda \rightarrow 0$, $\lambda \in \mathbb{C}' $:
\begin{equation} \label{Modd}
R^\beta_\lambda=Q^\beta_{d}
\lambda^{\frac{d}{2}-1}+{P^\beta_{d}}(\lambda)+O(|\lambda
|^{\frac{d}{2}}),~~~d \geq3,~~~d-{\rm odd},
\end{equation}
\begin{equation} \label{Meven}
R^\beta_\lambda=Q^\beta_{d} \lambda^{\frac{d}{2}-1} \ln
({1}/{\lambda}) +{P^\beta_{d}}(\lambda)+O(|\lambda
|^{\frac{d}{2}}|\ln \lambda|+ |\lambda|^{d-2}|\ln \lambda|^2
),~~~d \geq 4,~~~d-{\rm even},
\end{equation}
where $P^\beta_{d}$ are polynomials with coefficients that are
bounded operators and $Q^\beta_{d}$ are bounded operators.
\end{lemma}
%

The limit $P^\beta_d(0) = \lim_{\lambda \rightarrow 0, \lambda \in
\mathbb{C}'} R^\beta_\lambda$ will be denoted by $R^\beta_0$. It
is an operator acting from $C_{\exp}( \mathbb{R}^d)$ to $ C(
\mathbb{R}^d)$.
\\

It will be shown in the Appendix that if $\beta
> \beta_{\rm cr}$, then the eigenvalue $\lambda_0(\beta)$ of
the operator $L^\beta$ is simple and the corresponding
eigenfunction does not change sign (and so can be assumed to be
positive). From this and Lemma~\ref{rbetaW} it follows that the
residue of $R^\beta_\lambda$ at $\lambda_0$ is the integral
operator with the kernel $\psi_\beta(x) \psi_\beta(y)$, where
$\psi_\beta$ is the positive eigenfunction normalized by the
condition $||\psi_\beta||_{L^2( \mathbb{R}^d)} = 1$. Note that
$\psi_\beta$ decays exponentially at infinity. More precisely, it
follows from (\ref{ha2}) that if we write $x$ as $(\theta, |x|)$
in polar coordinates, then there is a continuous function
$f_\beta$ such that
\begin{equation} \label{uuo}
\psi_\beta(x) \sim f_\beta(\theta) |x|^{\frac{1}{2}-\frac{d}{2}}
\exp(-\sqrt{2\lambda_0}|x|)~~~{\rm as}~~|x| \rightarrow \infty.
\end{equation}
If $\beta = \beta_{\rm cr}$, then $\lambda_0 = 0$ might not be an
eigenvalue of the operator $L^\beta$. As shown in Lemma~\ref{re1},
for $d \geq 3$ there is a unique (up to a multiplicative constant)
positive function $\psi_{\beta}$ (the ground state of $L^\beta$)
which satisfies
\[
L^{\beta }\psi_{\beta} =\frac{1}{2}\Delta \psi +\beta
v\psi_{\beta} =0,~~~\psi_{\beta} (x) =
O(|x|^{2-d}),~~\frac{\partial{\psi_{\beta}} }{\partial{r}}(x) =
O(|x|^{1-d})~~as~~r =|x| \rightarrow \infty. \label{grsA}
\]
In fact, $\psi_\beta$ is a genuine eigenvector (element of
$L^{2}(\mathbb{R}^d)$) if and only if $d \geq 5$. We will
normalize $\psi_\beta$ by the condition that $||\beta v
\psi_\beta||_{L^2( \mathbb{R}^d)} = 1$.

\section{The super-critical case} \label{scc}

Throughout this section we assume that $\beta
> \beta_{\rm cr}$. First, let us introduce some notations.
Note that the dependence of some of the quantities below on
$\beta$ in not reflected in the notation in order to avoid
overcrowded formulas. For a positive number $x$, we define the
curve $\Gamma(x)$ in the complex plane as follows:
\[
\Gamma(x) = \{\lambda: |{\rm Im} \lambda| = \sqrt{ 4x(x -{\rm Re}
\lambda)},~{\rm Re} \lambda \geq 0 \} \bigcup \{\lambda: |{\rm Im}
\lambda| = 2x(1 -{\rm Re} \lambda),~{\rm Re} \lambda \leq 0 \}.
\]
Thus $\Gamma(x)$ is a union of a piece of the parabola with the
vertex in $x$ that points in the direction of the negative real
axis and two rays tangent to the parabola at the points it
intersects the imaginary axis. The choice of the curve is somewhat
arbitrary, yet the following properties of $\Gamma(x)$ will be
important:

First, ${\rm Re} \lambda \leq x$ for $\lambda \in \Gamma(x)$.
Second, since the rays form a positive angle with the negative
real semi-axis, we have $|{\rm arg} \lambda| \leq \pi -
\varepsilon (x)$ for all $\lambda \in \Gamma(x)$ for some
$\varepsilon(x) > 0$. Third, since the rays are tangent to the
parabola, and the parabola is mapped into the line $\{ \lambda:
{\rm Re} \lambda = \sqrt{x} \}$ by the mapping $\lambda
\rightarrow \sqrt{\lambda}$, the image of the curve $\Gamma(x)$
under the same mapping lies in the half-plane $\{\lambda: {\rm Re}
\lambda \geq \sqrt{x} \}$.

The integration along the vertical lines in the complex plane and
along contours $\Gamma(x)$, below, is performed in the direction
of the increasing complex part.

We'll need estimates on the solutions of the following parabolic
equation. Let
\begin{equation} \label{zeroind}
\partial_t \rho(t,x) = \frac{1}{2} \Delta \rho(t,x) +
\beta v(x) \rho(t,x),~~~~ \rho(0,x) = g(x) \in C_{\rm exp}.
\end{equation}
 We'll denote the Laplace transform of a function $f$ by $
\widetilde{f}$,
\[
\widetilde{f}(\lambda) = (\mathcal{L}f) (\lambda) =  \int_0^\infty
\exp(-\lambda t) f(t) d t.
\]
 Let $r$ be the distance between $\lambda_0$ and the rest of
the spectrum of the operator $L^\beta$. In the arguments that
follow we'll use the symbol $A$ to denote constants that may
differ from line to line.

\begin{lemma} \label{tth}
For each $\varepsilon \in(0, r)$, the solution of (\ref{zeroind})
has the form
\begin{equation} \label{tthh}
\rho(t,x) = \exp(\lambda_0 t) \langle \psi_\beta,  g \rangle
\psi_\beta(x) +   q(t,x),
\end{equation}
where
\[
||q(t,\cdot)||_{C}  \leq A(\varepsilon)
\exp((\lambda_0-\varepsilon) t) ||g||_{C_{\rm exp}}.
\]
\end{lemma}
\proof After the Laplace transform, the equation becomes
\[
(\frac{1}{2} \Delta + \beta v) \widetilde{\rho}  - \lambda
\widetilde{\rho} = - g.
\]
Thus, the solution $\rho$ can be represented as
\begin{equation} \label{r11}
\rho(t,\cdot) = -\frac{1}{2 \pi i} \int_{{{\rm Re} \lambda =
\lambda_0 +1}} e^{\lambda t} R^\beta_\lambda  g d \lambda.
\end{equation}
The resolvent is meromorphic in the complex plane outside of the
interval $(-\infty, \lambda_0-r]$,  with the only (simple) pole
at~$\lambda_0$
 with the principal part of the Laurent expansion being
the integral operator with the kernel $\psi_\beta(x)
\psi_\beta(y)/(\lambda_0 - \lambda)$.

By Lemma~\ref{rbetaW}, the norm of the $R^\beta_\lambda$ does not
exceed $A/|\lambda|$ near infinity to the right of
$\Gamma(\lambda_0 -\varepsilon)$. Therefore, the same integral as
in (\ref{r11}) but along the segment parallel to the real axis
connecting a point $\lambda_0 + 1 + i b$ with the contour
$\Gamma(\lambda_0 -\varepsilon)$   tends to zero when $b
\rightarrow \infty$. Therefore, we can replace the contour of
integration in (\ref{r11}) by $\Gamma(\lambda_0 -\varepsilon)$.
The residue gives the main term, while
the integral over $\Gamma(\lambda_0 -\varepsilon)$ gives the
remainder term. \qed

\begin{lemma} \label{rhoo} Let $K \subset \mathbb{R}^d$ be a compact set. For
each $\varepsilon \in (0,r)$, the function $\rho_1(t,x,y)$
satisfies
\[
\rho_1(t,x,y) = \exp( \lambda_0 t) \psi_\beta(x) \psi_\beta(y)  +
q (t,x,y),
\]
where
\[
\sup_{x \in K} |q(t,x,y)| \leq A(\varepsilon) \exp((\lambda_0-
\varepsilon) t - |y| \sqrt{2 (\lambda_0 - \varepsilon)})
\]
for $t \geq 1/2$.
\end{lemma}
\proof Let $K'$ be a compact set that contains ${\rm supp}(v) \cup
K $ in its interior. Consider first the case when $y \in K'$.
Apply (\ref{tthh}) with $t$ replaced by $t' = t - 1/2$ and $g =
\rho_1(1/2,\cdot,y)$. In order to calculate the main term of the
asymptotics, we note that $ ||g||_{C_{\rm exp}}$ is bounded
uniformly in $y \in K'$ and
\[
\langle \psi_\beta, g \rangle = \exp(\frac{1}{2} {\lambda_0})
\psi_\beta(y).
\]
The latter follows from
\[
0 = \int_0^{1/2} \langle (\frac{\partial}{\partial t} +
L^\beta)(\exp(-\lambda_0 t) \psi_\beta), \rho_1 \rangle d t
=
\]
\[
\langle \exp(-\lambda_0 t) \psi_\beta, \rho_1 \rangle|_{t=0}^{1/2}
+ \int_0^{1/2} \langle (\exp(-\lambda_0 t) \psi_\beta), (-
\frac{\partial}{\partial t} + L^\beta) \rho_1 \rangle d t =
\]
\[
\langle \exp(-\frac{1}{2}\lambda_0) \psi_\beta,
\rho_1(1/2,\cdot,y) \rangle - \langle  \psi_\beta, \rho_1
(0,\cdot,y) \rangle =
\]
\[
\exp(-\frac{1}{2}\lambda_0)  \langle \psi_\beta, g \rangle -
\psi_\beta(y).
\]
Therefore, (\ref{tthh}) implies that
\[
\rho_1(t,x,y) = \exp(\lambda_0 t)\psi_\beta(y) \psi_\beta(x) +
\exp((\lambda_0 -\varepsilon) t ) q(t,x,y),
\]
where $||q(t,\cdot,y)||_C \leq A(K')$ for all $y \in K'$. It
remains to consider the case when $y \notin K'$.

Let $u(t,x,y) = \rho_1(t,x,y) - p_0(t,x,y)$, where $p_0$ is the
fundamental solution of the heat equation. Then $u$ satisfies the
non-homogeneous version of (\ref{zeroind}) with the right hand
side $f = -\beta v(x) p_0(t,x,y)$ and $g \equiv 0$. Note that $f$
is a smooth function since $y \notin K'$. Solving this equation
for $u$ using the Laplace transform, as in the proof of
Lemma~\ref{tth}, we obtain
\[
u(t,\cdot,y) = -\frac{1}{2 \pi i} \int_{{{\rm Re} \lambda =
\lambda_0 +1}} e^{\lambda t} R^\beta_\lambda  (-\beta v
\widetilde{p}_0(\lambda, \cdot, y)) d \lambda
\]
\begin{equation} \label{iioo}
= - \frac{1}{2 \pi i} \int_{{{\rm Re} \lambda = \lambda_0 +1}}
e^{\lambda t} R^\beta_\lambda  (\beta v R^0_{\lambda}(\cdot,y)) d
\lambda
\end{equation}
\[
=\exp(\lambda_0 t)\langle \psi_\beta, \beta v
R^0_{\lambda_0}(\cdot,y) \rangle \psi_\beta - \frac{1}{2 \pi i}
\int_{ \Gamma( \lambda_0 -\varepsilon)} e^{\lambda t}
R^\beta_\lambda  (\beta v R^0_{\lambda}(\cdot,y)) d \lambda,
\]
where the first term on the right hand side is due to the residue
at $\lambda = \lambda_0$. The first term can be re-written as
\[
\exp(\lambda_0 t)\langle \psi_\beta, \beta v
R^0_{\lambda_0}(\cdot,y) \rangle \psi_\beta(x) =
\]
\[
\exp(\lambda_0 t)(R^0_{\lambda_0} (\beta v \psi_\beta)) (y)
\psi_\beta(x)  = -\exp(\lambda_0 t)\psi_\beta(y) \psi_\beta(x).
\]
The last equality here follows from the fact that $\psi_\beta$ is
an eigenfunction with eigenvalue $\lambda_0$, that is
\[
(\frac{1}{2} \Delta  - \lambda_0) \psi_\beta = -\beta v
\psi_\beta.
\]
In order to estimate the second term on the right hand side of
(\ref{iioo}), we note that from~(\ref{ha2}) (see also
(\ref{estwT})) it follows that
\[
|R^{0}_\lambda(x,y)| \leq  A(l)
|\sqrt{\lambda}|^{\frac{d}{2}-\frac{3}{2}} |x - y
|^{\frac{1}{2}-\frac{d}{2}}  |\exp(-\sqrt{2 \lambda} |y-x|)|
\]
if $|\lambda|, |y -x| \geq l$. Thus
\[
 ||\beta v
R^0_{\lambda}(\cdot,y)||_{C_{\rm exp}} \leq A(\varepsilon) |y
|^{\frac{1}{2}-\frac{d}{2}}
|\sqrt{\lambda}|^{\frac{d}{2}-\frac{3}{2}} \exp(-\sqrt{2(\lambda_0
- \varepsilon)} |y|)
\]
for $y \notin K'$, $\lambda \in \Gamma(\lambda_0 - \varepsilon)$
due to the fact that ${\rm Re} \sqrt{\lambda} \geq \sqrt{\lambda_0
- \varepsilon}$ for $\lambda \in \Gamma(\lambda_0 - \varepsilon)$
and $|y - x| \geq l$ for $x \in {\rm supp}(v)$, $y \notin K'$.

Hence, using the estimate on the norm of $R^\beta_\lambda : C_{\rm
exp} \rightarrow C$ from Lemma~\ref{rbetaW}, we obtain
\[
|| R^\beta_\lambda  (\beta v R^0_{\lambda}(\cdot,y)) ||_C \leq
A(\varepsilon) |\sqrt{\lambda}|^{\frac{d}{2}-\frac{5}{2}}
\exp(-\sqrt{2(\lambda_0 - \varepsilon)} |y|),~~~~\lambda \in
\Gamma(\lambda_0 - \varepsilon).
\]
Therefore, since ${\rm Re} \lambda \leq \lambda_0 -\varepsilon$
for $\lambda \in \Gamma(\lambda_0 -\varepsilon)$ and the factor
$e^{\lambda t}$ decays exponentially along $\Gamma(\lambda_0
-\varepsilon)$, the $C$-norm of the second term on the right hand
side of (\ref{iioo}) does not exceed $A(\varepsilon)
\exp((\lambda_0- \varepsilon) t - |y| \sqrt{2 (\lambda_0 -
\varepsilon)})$. The term $p_0(t,x,y)$ with $x \in K$, $y \notin
K'$, $t \geq 1/2$, is estimated by the same expression, possibly
with a different constant $A(\varepsilon )$. Indeed, if $t \geq
1/2$, then
\[
p_0(t,x,y) \leq A \exp(-|y-x|^2/2t) \leq A \exp((\lambda_0-
\varepsilon) t - |y-x| \sqrt{2 (\lambda_0 - \varepsilon)})
\]
since
\[
|y-x|^2/2t + (\lambda_0- \varepsilon) t - |y-x| \sqrt{2 (\lambda_0
- \varepsilon)} = (|y-x|/\sqrt{2t} - \sqrt{ (\lambda_0-
\varepsilon) t} ) ^2 \geq 0.
\]
\qed

We'll  need additional notations in order to describe the
asymptotics of $\rho_n$ with $n > 1$. Let
$\alpha^1_\varepsilon(t,y) = \psi_\beta(y)$ and
$\alpha^2_\varepsilon(t,y) = \exp(-\varepsilon t - |y| \sqrt{2
(\lambda_0 - \varepsilon)})$. Consider all possible sequences
$\sigma = (\sigma_1,...,\sigma_n)$ with $\sigma_i \in \{1,2\}$. By
$\Pi^n_\varepsilon (t,y_1,...,y_n)$ we denote the quantity
\[
\Pi^n_\varepsilon (t,y_1,...,y_n) = \sup_{\sigma \neq (1,...,1)}
\alpha^{\sigma_1}_\varepsilon(t,y_1)\cdot ... \cdot
\alpha^{\sigma_n}_\varepsilon(t,y_n).
\]

Let $P_t: C_{\rm \exp} \rightarrow C$ be the operator that maps
the initial function $g$ to the solution $\rho(t,\cdot)$ of
equation (\ref{zeroind}). Let $P^0_t g (x) = \exp(\lambda_0 t)
\langle \psi_\beta,  g \rangle \psi_\beta(x)$ and $P^1_t = P_t -
P^0_t$. Lemma~\ref{tth} states that
\[
||P^1_t|| \leq A(\varepsilon)
\exp((\lambda_0-\varepsilon) t).
\]
The particular form of $P^0_t$ then implies that
\begin{equation} \label{hjh}
||P_t|| \leq  ||P^0_t|| + ||P^1_t|| \leq A'  \exp(\lambda_0 t).
\end{equation}
For $g \in C_{\rm exp}$ and $n \geq 2$,  we denote
\[
I_n(g) := R^\beta_{ n \lambda_0} g  = \int_0^\infty \exp(-n
\lambda_0 s) P_s g d s \in C.
\]
Note that
\[
\int_0^t  \exp(n \lambda_0 s) P_{t-s} g d s = \exp(n \lambda_0 t)
\int_0^t  \exp(-n \lambda_0 s) P_s g d s
\]
\begin{equation} \label{amm}
= \exp(n \lambda_0 t)(I_n(g) + O(\exp(-(n-1) \lambda_0 t))
)~~~~{\rm as}~~t \rightarrow \infty.
\end{equation}
The functions $f_1, f_2,...$  are defined inductively: $f_1 =
\psi_\beta$ and
\begin{equation} \label{fnh}
f_n =  \beta \sum_{k=1}^{n-1} \frac{n!}{k! (n-k)!} I_n(v f_k
f_{n-k}),~~~n \geq 2.
\end{equation}

\begin{lemma} \label{rnn}
Let $K \subset \mathbb{R}^d$ be a compact set. For each
$\varepsilon \in (0,r)$, the function $\rho_n$ satisfies
\begin{equation} \label{yuy2}
\rho_n(t,x,y_1,...,y_n) =  \exp( n \lambda_0 t) f_n(x)
\psi_\beta(y_1) \cdot ... \cdot \psi_\beta(y_n) +
q_n(t,x,y_1,...,y_n),
\end{equation}
where
\begin{equation} \label{ggp}
\sup_{x \in K} |q_n(t,x,y_1,...,y_n)| \leq A_n(\varepsilon) \exp(n
\lambda_0 t ) \Pi^n_\varepsilon (t,y_1,...,y_n)
\end{equation}
for $t \geq 1/2$.
\end{lemma}
\proof For $n =1$, the relation (\ref{yuy2}) coincides with the
statement of Lemma~\ref{rhoo}. Let us assume that (\ref{yuy2})
holds for all natural numbers up to and including $n-1$. A generic
subsequence $U \subset Y = (y_1,...,y_n)$ will be written as $U =
(z_1,...,z_{|U|})$ and its complement as $Y \setminus U =
(\overline{z}_1,...,\overline{z}_{n - |U|})$. By the Duhamel
principle applied to the equation for $\rho_n$, we obtain
\[
\rho_n(t,\cdot,y_1,...,y_n) = \int_0^t P_{t-s} (\beta v \sum_{U
\subset Y, U \neq \emptyset} \rho_{|U|} (s,\cdot,z_1,...,z_{|U|})
\rho_{n-|U|}(s,\cdot,\overline{z}_1,...,\overline{z}_{n - |U|})) d
s
\]
\[
= \int_0^t P_{t-s} (\beta v  \sum_{U \subset Y, U \neq \emptyset}
\exp( |U| \lambda_0 s) f_{|U|} (\cdot)
\psi_\beta(z_1)\cdot...\cdot \psi_\beta(z_{|U|}))
\]
\begin{equation} \label{middxx}
\times  \exp( (n-|U|) \lambda_0 s) f_{n - |U|} (\cdot)
\psi_\beta(\overline{z}_1)\cdot...\cdot
\psi_\beta(\overline{z}_{n-|U|})) d s
\end{equation}
\[
+  2 \int_0^t P_{t-s} (\beta v  \sum_{U \subset Y, U \neq
\emptyset} \exp( |U| \lambda_0 s) f_{|U|} (\cdot)
\psi_\beta(z_1)\cdot...\cdot \psi_\beta(z_{|U|})
q_{n-|U|}(s,\cdot,\overline{z}_1,...,\overline{z}_{n - |U|})) d s
\]
\[
+ \int_0^t P_{t-s} (\beta v  \sum_{U \subset Y, U \neq \emptyset}
q_{|U|}(s,\cdot,{z}_1,...,{z}_{|U|})
q_{n-|U|}(s,\cdot,\overline{z}_1,...,\overline{z}_{n - |U|})) d s.
\]
The second and third integrals on the right hand side of
(\ref{middxx}) contribute only to the remainder term. Indeed,
consider the contribution to the second integral from the term
with a given $U$:
\[
\int_0^t P_{t-s} (\beta v  \exp( |U| \lambda_0 s) f_{|U|} (\cdot)
\psi_\beta(z_1)\cdot...\cdot \psi_\beta(z_{|U|})
q_{n-|U|}(s,\cdot,\overline{z}_1,...,\overline{z}_{n - |U|})) d s
\]
\[
\leq A \psi_\beta(z_1)... \psi_\beta(z_{|U|}) \int_0^t P_{t-s} ( v
\exp( |U| \lambda_0 s) f_{|U|} (\cdot)
\]
\[
\times
 \exp((n - |U|)
\lambda_0 s ) \Pi^{n-|U|}_\varepsilon
(s,\overline{z}_1,...,\overline{z}_{n - |U|})) d s
\]
\[
\leq  A \psi_\beta(z_1)\cdot...\cdot \psi_\beta(z_{|U|}) \int_0^t
\exp(\lambda_0(t-s)) \exp(n \lambda_0 s) \Pi^{n-|U|}_\varepsilon
(s,\overline{z}_1,...,\overline{z}_{n - |U|})) d s
\]
\[
\leq  A \exp(n \lambda_0 t) \psi_\beta(z_1)... \psi_\beta(z_{|U|})
\Pi^{n-|U|}_\varepsilon (t,\overline{z}_1,...,\overline{z}_{n -
|U|}) \leq  A\exp(n \lambda_0 t ) \Pi^n_\varepsilon
(t,y_1,...,y_n),
\]
where the first inequality follows from the inductive assumption
and the second one from~(\ref{hjh}). The third integral on the
right hand side of (\ref{middxx}) is estimated similarly. It
remains to consider the first integral. It is equal to
\[
\psi_\beta(y_1)\cdot...\cdot \psi_\beta(y_n)
 \int_0^t   \exp (n \lambda_0 s) P_{t-s} (\beta
v    \sum_{U \subset Y, U \neq \emptyset} f_{|U|} f_{n-|U|} ) d s
\]
\[
= \psi_\beta(y_1) \cdot ... \cdot \psi_\beta(y_n)   \exp( n
\lambda_0 t) \left( f_n(\cdot) +  O(\exp(-(n-1) \lambda_0
t))\right),
\]
where the last equality follows from (\ref{amm}). Thus we obtain
the main term from the right hand side of (\ref{yuy2}) plus the
correction
\[
\psi_\beta(y_1) \cdot ... \cdot \psi_\beta(y_n)   \exp( n
\lambda_0 t)   O(\exp(-(n-1) \lambda_0 t))
\]
for which the estimate (\ref{ggp}) holds since $ \psi_\beta(y_1)
\exp(-\lambda_0 t) \leq  \alpha^2_\varepsilon(t,y_1)$ due to
(\ref{uuo}). \qed
\\

We will now  use Lemma~\ref{rnn} to draw conclusions about the
distribution of particles at time $t$. First, let us observe that
for each $x$ the sequence $\left( \int_{ \mathbb{R}^d}
\psi_\beta(y) d y \right)^n f_n(x)$, $n \geq 1$,  serves as a
sequence of moments for a random variable $\xi^{\beta,x}$ whose
distribution is defined uniquely. Indeed, by the Carleman theorem,
it is sufficient to check that
\begin{equation} \label{fnhx}
\sum_{n =1}^\infty \left( \frac{1}{f_n(x)} \right)^{\frac{1}{2n}}
= \infty.
\end{equation}
From (\ref{hjh}) it follows that there is a constant $A$ such that
\[
|| I_n(g) ||_C \leq \frac{A}{n-1} ||g||_{C_{\rm exp}},~~~n \geq 2.
\]
Therefore, from (\ref{fnh}) it follows that for a different
constant $A$,
\[
||f_n||_C \leq  A \sum_{k=1}^{n-1} \frac{(n-1)!}{k! (n-k)!}
||f_k||_C ||f_{n-k}||_C,~~~n \geq 2,~~~~~||f_1||_C \leq A.
\]
From here, by induction on $n$ it follows that $||f_n||_C \leq
A^{2n-1} n!$, which in turn implies (\ref{fnhx}) since $n! \leq
((n+1)/2)^n$.

Let $n^{\beta,x}_t(U)$ be the number of particles in a domain $U
\subseteq \mathbb{R}^d$, assuming that at $t = 0$ there was one
particle located at $x$. We will write $n^{\beta,x}_t$ instead of
$n^{\beta,x}_t( \mathbb{R}^d)$. Note that
\begin{equation} \label{stf}
{\mathrm{E} }( n^{\beta,x}_t(U) )^n = \sum_{k = 1}^n S(n,k)
\int_U...\int_U \rho_k(t,x,y_1,...,y_k) dy_1...dy_k,
\end{equation}
where $S(n,k)$ is the Stirling number of the second kind (the
number of ways to partition $n$ elements into $k$ nonempty
subsets). Formula (\ref{stf}) easily follows if we write
\[
n^{\beta,x}_t(U) = \sum_i n^{\beta,x}_t(\Delta_i),
\]
where $U = \sqcup_i \Delta_i$ is the partition of $U$ into small
sub-domains, and then take the limit as $\max_i {\rm diam}
(\Delta_i) \rightarrow 0$.

Let $\xi^{\beta,x}$ be a random variable with the moments
\[
{\mathrm{E}}(\xi^{\beta,x})^n = \left( \int_{ \mathbb{R}^d}
\psi_\beta(y) d y \right)^n f_n(x).
\]
Let $\varphi_\beta$ be the density on $ \mathbb{R}^d$ given by
\[
\varphi_\beta(x) = \psi_\beta(x)/\int_{ \mathbb{R}^d}
\psi_\beta(y) d y.
\]
\begin{theorem} \label{mtj}
For each $x \in \mathbb{R}^d$ and each domain $U \subseteq
\mathbb{R}^d$,
\[
\lim_{t \rightarrow \infty} \left( \exp(-\lambda_0 t)
n^{\beta,x}_t(U) \right) = \xi^{\beta,x} \int_U \varphi_\beta(x) d
x
\]
in distribution.
\end{theorem}
\proof By the theorem of Frechet and Shohat \cite{FSh}, it is
sufficient to prove the convergence of the moments. By (\ref{stf})
the $n$-th moment of $\exp(-\lambda_0 t) n^{\beta,x}_t(U)$ is
equal to
\[
{\mathrm{E} }(\exp(-\lambda_0 t)n^{\beta,x}_t(U))^n = \exp(-n
\lambda_0 t) \sum_{k = 1}^n S(n,k) \int_U...\int_U
\rho_k(t,x,y_1,...,y_k) dy_1...dy_k.
\]
First consider the contribution to the right hand side from the
term with $k = n$. Note that $S(n,n) = 1$.  We use (\ref{yuy2})
for the asymptotics of $\rho_n$. The contribution from the term
$q_n(t,x,y_1,...,y_n)$ tends to zero:
\[
\lim_{t \rightarrow \infty} \exp(-n \lambda_0 t)\int_U...\int_U
q_n(t,x,y_1,...,y_n) dy_1...dy_n = 0,
\]
as follows from the definition of $\Pi^n_\varepsilon
(t,y_1,...,y_n)$. The contribution from the main term gives the
desired expression:
\[
\lim_{t \rightarrow \infty} \exp(-n \lambda_0 t)\int_U..\int_U
\exp( n \lambda_0 t) f_n(x) \psi_\beta(y_1)  ... \psi_\beta(y_n)
dy_1...dy_n = \mathrm{E}(\xi^{\beta,x})^n ( \int_U
\varphi_\beta(x) d x )^n.
\]
It remains to note that the contribution from each of the terms
with $k < n$ tends to zero. Indeed, it is equal to
\[
\exp(-(n-k) \lambda_0 t) \left(   \exp(-k \lambda_0 t) S(n,k)
\int_U...\int_U \rho_k(t,x,y_1,...,y_k) dy_1...dy_k \right).
\]
The expression inside the brackets tends to a finite limit as in
the case $k = n$, while the exponential factor in front of the
brackets tends to zero for $k < n$.
 \qed

\section{The sub-critical case}

Throughout this section we assume that $d \geq 3$ and $\beta \in
(0, \beta_{\rm cr})$. We will show that  the total number of
particles $n^{\beta,x}_t$ tends to a random limit when $t
\rightarrow \infty$. Denote the integrals of the correlation
functions by $\overline{\rho}_n(t,x)$:
\[
\overline{\rho}_n(t,x) =
\int_{\mathbb{R}^d}...\int_{\mathbb{R}^d}\rho_n(t,x,y_1,...,y_n)
dy_1....dy_n.
\]
From (\ref{fmo})-(\ref{ico}) and (\ref{manyp})-(\ref{initu}) it
follows that these quantities satisfy the equations
\[
\partial_t \overline{\rho}_1(t,x) = \frac{1}{2} \Delta \overline{\rho}_1(t,x) +
\beta v(x) \overline{\rho}_1(t,x),
\]
\[
\overline{\rho}_1(0,x) \equiv  1
\]
and for $n > 1$:
\[
\partial_t \overline{\rho}_n(t,x) = \frac{1}{2} \Delta \overline{\rho}_n(t,x) +
\beta v(x) \left( \overline{\rho}_n(t,x) + \sum_{k =1}^{n-1}
\frac{n!}{k! (n-k)!} \overline{\rho}_k (t,x) \overline{\rho}_{n-k}
(t,x) \right) ,
\]
\[
\overline{\rho}_n(0,x) \equiv 0.
\]
As in the previous section, in order to find the asymptotics of
the total number of particles, we will study the asymptotics of
$\overline{\rho}_n(t,x)$ as $r \rightarrow \infty$.

\begin{lemma} \label{tth2}
For $\beta \in (0, \beta_{\rm cr})$,  there is a constant $A$ such
that the solution $\rho$ of (\ref{zeroind}) can be estimated as
follows
\begin{equation} \label{nnj}
||\rho(t,\cdot)||_C \leq A (1+t)^{-d/2} ||g||_{C_{\rm exp}}.
\end{equation}
If the initial condition $g \in C_{\rm exp}$ in (\ref{zeroind}) is
replaced by $g \equiv 1$, then
\begin{equation} \label{zic}
\lim_{t \rightarrow \infty} \rho(t,x) =  1+ \varphi_\beta(x)
\end{equation}
in $C(\mathbb{R}^d)$, where $\varphi_\beta = - R^\beta_0(\beta
v)$.
\end{lemma}
\proof Let $d \geq 3$ be odd. As in (\ref{r11}), we represent the
solution as an integral
\[
\rho(t,\cdot) = -\frac{1}{2 \pi i} \int_{{{\rm Re} \lambda = 1}}
e^{\lambda t} R^\beta_\lambda  g d \lambda.
\]
Since the integrand is analytic in $ \mathbb{C}'$ with appropriate
decay at infinity and has a limit as $\lambda \rightarrow 0$, we
can replace the contour of integration by
\[
\gamma = \{ z \in \mathbb{C}: {\rm Re} z = -|{\rm Im} z| \}.
\]
Using the representation (\ref{Modd}) for $R^\beta_\lambda$, we
obtain
\[
\rho(t,\cdot) = -\frac{1}{2 \pi i} \int_\gamma e^{\lambda t}(
Q^\beta_{d}
\lambda^{\frac{d}{2}-1}+{P^\beta_{d}}(\lambda)+O(|\lambda
|^{\frac{d}{2}}) ) g d \lambda.
\]
The contribution to the integral from the first term is estimated
as follows
\[
||Q^\beta_{d} g \int_\gamma e^{\lambda t} \lambda^{\frac{d}{2}-1}
d \lambda ||_C =  || Q^\beta_{d} gt^{-\frac{d}{2}} \int_\gamma
e^{s} s^{\frac{d}{2}-1} d s ||_C \leq A t^{-\frac{d}{2}}
||g||_{C{\rm exp}},
\]
where we used the change of variable $s  = \lambda t$. The
contribution from the second term is equal to zero since the
integrand is analytic and the contour can be moved arbitrarily far
to the left along the real axis. The third term can be treated in
the same way as the first one, resulting in the decay in $t$ of
order $ t^{-\frac{d}{2}-1}$. The obtained estimates imply
(\ref{nnj}) for $t \geq 1$. Clearly (\ref{nnj}) holds for $t \leq
1$.

The case when $d \geq 4$ is even is treated similarly. The slight
difference is that the contribution from the main term is now, up
to a multiplicative constant, equal to
\[
Q^\beta_{d} g \int_\gamma e^{\lambda t} \lambda^{\frac{d}{2}-1}
\ln(1/\lambda) d \lambda.
\]
After the change of variable $s  = \lambda t$, the integral is
seen to be equal to
\[
t^{-\frac{d}{2}} \int_\gamma e^{s} s^{\frac{d}{2}-1} \ln(t/s) d s
= t^{-\frac{d}{2}} \ln t \int_\gamma e^{s} s^{\frac{d}{2}-1} d s +
t^{-\frac{d}{2}} \int_\gamma e^{s} s^{\frac{d}{2}-1} \ln(1/s) d s.
\]
The first integral on the right hand side is equal to zero since
the integrand is an analytic function, and the contour can
therefore be moved arbitrarily far to the left along the real
axis. The second term on the right hand side has the desired order
in $t$.

It remains to prove (\ref{zic}). Note that $w(t,x) = \rho(t,x)- 1$
is the solution of the problem
\[
\frac{\partial w(t,x)}{\partial t} = \frac{1}{2}\Delta w(t,x) +
\beta v(x) w(t,x) + \beta v(x),~~~w(0,x) \equiv 0.
\]
By the Duhamel formula,
\[
 w(t,\cdot) = \frac{-1}{2 \pi i} \int_0^t \int_{\gamma} e^{\lambda( t-s)}
R^\beta_\lambda (\beta v) d \lambda d s =
\]
\[
\frac{-1}{2 \pi i}\int_{\gamma} \frac{e^{\lambda t}-1}{\lambda}
R^\beta_\lambda( \beta v) d \lambda = \frac{-1}{2 \pi
i}\int_{\gamma} \frac{e^{\lambda t}}{\lambda} R^\beta_\lambda
(\beta v) d \lambda,
\]
since in the domain $\gamma^+$ to the right of the contour
$\gamma$, the operator $ R^\beta_\lambda$ is analytic and decays
as $|\lambda|^{-1}$ at infinity. We make the change of variables
$s = \lambda t$ and obtain, as $t \rightarrow \infty$,
\[
\frac{-1}{2 \pi i}\int_{\gamma} \frac{e^{\lambda t}}{\lambda}
R^\beta_\lambda (\beta v) d \lambda = \frac{-1}{2 \pi
i}\int_{\gamma} \frac{e^{s}}{s} R^\beta_{s/t} (\beta v) d s
\rightarrow R^\beta_{0} (\beta v) \frac{-1}{2 \pi i}\int_{\gamma}
\frac{e^{s}}{s}  d s = \varphi_\beta.
\]
 \qed
\\

For $g \in C_{\rm exp}$, we define $ J(g) = \int_0^\infty P_s g d
s$. From Lemma~\ref{tth2} it follows that
\begin{equation} \label{twee}
||J(g)||_C \leq A ||g||_{C_{\rm exp}}
\end{equation}
for some constant $A$. Let us define the sequence of functions
$f_1, f_2,...$ inductively via: $f_1 = 1 + \varphi_\beta$ and
\begin{equation} \label{fnh2}
f_n =   \beta \sum_{k=1}^{n-1} \frac{n!}{k! (n-k)!}  J( v f_k
f_{n-k}),~~~n \geq 2.
\end{equation}
\begin{lemma} \label{cmo}
For each $n \geq 1$ we have
\begin{equation} \label{limi2}
\lim_{t \rightarrow \infty} \overline{\rho}_n(t,\cdot)  = f_n
\end{equation}
in $ C( \mathbb{R}^d)$.
\end{lemma}
\proof For $n =1$, the statement coincides with the second part of
of Lemma~\ref{tth2}. Let us assume that (\ref{limi2}) holds for
all natural numbers up to and including $n-1$. Define the
functions $c_k(t,x) = \overline{\rho}_k(t,x) - f_k(x)$. Thus
$||c_k(t,\cdot)||_C \rightarrow 0$ as $t \rightarrow \infty$ for
$k \leq n-1$. By the Duhamel principle applied to the equation for
$\overline{\rho}_n$, we obtain
\[
\overline{\rho}_n(t,\cdot) = \int_0^t P_{t-s} (\beta v \sum_{k
=1}^{n-1}  \frac{n!}{k! (n-k)!}  \overline{\rho}_k (s,\cdot)
\overline{\rho}_{n-k}(s,\cdot) )d s
\]
\begin{equation} \label{middxx2}
=  \int_0^t P_{t-s} (\beta v \sum_{k =1}^{n-1}   \frac{n!}{k!
(n-k)!}  f_k  f_{n-k})d s +
\end{equation}
\[
 2 \int_0^t P_{t-s} (\beta v \sum_{k =1}^{n-1}   \frac{n!}{k!
(n-k)!}  c_k (s,\cdot) f_{n-k})d s + \int_0^t P_{t-s} (\beta v
\sum_{k =1}^{n-1} \frac{n!}{k! (n-k)!}  c_k (s,\cdot)
c_{n-k}(s,\cdot) )d s.
\]
From Lemma~\ref{tth2} it follows that $||P_t g||_C \leq A
(1+t)^{-d/2} ||g||_{C_{\rm exp}}$ for some constant $A$. Therefore
the $C$-norm of the sum of the last two term on the right hand
side of (\ref{middxx2}) is estimated from above by
\[
\int_0^t \gamma(s) (1+ t-s)^{-d/2} d s,
\]
where $\gamma(s)$ is some function such that $\gamma(s)
\rightarrow 0$ as $s \rightarrow \infty$. The latter integral
tends to zero when $t \rightarrow \infty$ since $d \geq 3$. It
remains to note that the first term on the right hand side of
(\ref{middxx2}) tends to $f_n$ in the $C$-norm, as immediately
follows from the definition of $J$. \qed
\begin{theorem} \label{mtk}
For each $x \in \mathbb{R}^d$, the total number of particles
$n^{\beta,x}_t$ converges almost surely, as $t \rightarrow
\infty$, to a random variable $\zeta^{\beta,x}$ with the moments
$m_n(x) = \sum_{k=1}^n S(n,k) f_k(x)$.
\end{theorem}
\proof Since $n^{\beta,x}_t$ is monotonically increasing in $t$,
and the moments converge to $m_n(x)$ (as follows from (\ref{stf})
and Lemma~\ref{cmo}), we have convergence almost surely to a
random variable with the specified moments. \qed

\section{The critical case}
\label{critc}
 Throughout this section we assume
that $d \geq 3$. We will show that  for $\beta = \beta_{\rm cr}$
the total number of particles tends almost surely to a finite
random limit, while the expectation of the number of particles
tends to infinity.

It is clear that the random variables $n^{\beta,x}_t$ can be
realized on a common probability space in such a way that
$n^{\beta,x}_t \leq n^{\beta',x}_{t'}$ whenever $\beta \leq
\beta'$ and $t \leq t'$. Therefore, in order to show that $
\mathrm{E} n^{\beta_{\rm cr},x}_t \rightarrow \infty$ as $t
\rightarrow \infty$, it is sufficient to show that
\[
\lim_{\beta \uparrow \beta_{\rm cr}} \lim_{t \rightarrow \infty}
\mathrm{E} n^{\beta,x}_t = \infty.
\]
This follows from Lemma 7.3 of \cite{CKMV1}, which can be stated
as follows:
\begin{lemma}
There are positive constant $b_d$, $d \geq 3$, such that
\[
\lim_{t \rightarrow \infty} \mathrm{E} n^{\beta,x}_t =
\frac{b_d}{\beta_{cr} - \beta} \psi_{\beta_{\rm cr}}(x) + O
(1)~~as~~\beta \uparrow \beta_{cr}
\]
is valid in $C( \mathbb{R}^d)$, where $\psi_{\beta_{\rm cr}}$ is
the positive ground state for $L^{\beta_{cr}}$.
\end{lemma}
Now we prove the main theorem of this section.
\begin{theorem}
For each $x$, the limit  $\eta^x := \lim_{t \rightarrow \infty}
n^{\beta_{\rm cr},x}_t $ is finite almost surely.
\end{theorem}
\proof Let $S(x) = \mathrm{P}(\eta^x = \infty)$. Let us show that
$S(x)$ depends continuously on $x$. Indeed, let $x \in
\mathbb{R}^d$ and $\varepsilon > 0$ be fixed. There exist $\delta
> 0$ and $t > 0$ such that for $|y-x| \leq \delta$ the
branching processes starting at $x$ and $y$ have the following
properties:

(a) The probability that at least one branching occurs on the
interval $[0,t]$ for either of the processes does not exceed
$\varepsilon/3$;

(b) There is a positive function $q(z)$ such that $\int_{
\mathbb{R}^d} q(z) \geq 1 - \varepsilon/3$, which serves as a
lower bound for both of the heat kernels $p(t,x,z)$ and
$p(t,y,z)$.

This shows that the branching processes starting at $x$ and $y$
can be coupled on an event of probability at least $1 -
\varepsilon$, and therefore $|S(x) - S(y)| \leq \varepsilon$,
proving the continuity.

Suppose that $S(x)$ is not identically equal to zero. Then $S(x)
\neq 0$ on a set of positive measure, due to the continuity. By
the Markov property, $S(x) \geq \int_{ \mathbb{R}^d} p(1,x,z) S(z)
d z$, which shows that $S(x) > 0$ for all $x \in  \mathbb{R}^d$.
Let
\[
m = \min_{x \in {\rm supp}(v)} S(x) > 0.
\]
Let $\Omega^x = \{ \eta^x = \infty \}$. Take $N = [4/m]+1$. Let
$r^x$ be the probability that there are at least $N$ particles on
the support of $v$ before time $t = 1$. Clearly, there is a
positive constant $r$ such that
\begin{equation} \label{mpro}
r^x \geq r~~~{\rm for}~~{\rm all}~~x \in {\rm supp} (v).
\end{equation}

Define the random times $\tau^x_1, \tau^x_2,...$ to be the
consecutive instances of branching for the process starting at
$x$. There is at least one particle on ${\rm supp} (v)$ at each of
the times $\tau^x_n$, and $\tau^x_n \rightarrow \infty$ almost
surely on $\Omega^x$. Therefore, from (\ref{mpro}) and the Markov
property it follows that  almost surely on $\Omega^x$ there is a
random time $\tau^x$ such that there are $N$ particles on the
support of $v$ at the time $\tau^x$. Therefore, there exists $T <
\infty$ such that $ \mathrm{P}(\tau^x \leq T) \geq m/2$ for $x \in
{\rm supp} (v)$. Note that $T$ can be taken to be independent of
$x$ using the continuity in~$x$ of the probabilities under
consideration and the compactness of ${\rm supp} (v)$.

Now fix an arbitrary $x \in {\rm supp} (v)$. We saw that with
probability at least $m/2$ there are at least $N$ particles at
time $\tau^x$ and therefore also at time $T$, that is $ \mathrm{E}
n^{\beta_{\rm cr},x}_T \geq Nm/2$. Applying the Markov property
with respect to the stopping time $\tau^x$ and using the fact that
the particles move independently, we see that $ \mathrm{E}
n^{\beta_{\rm cr},x}_{2T} \geq (Nm/2)^2$ and, continuing by
induction, that $ \mathrm{E} n^{\beta_{\rm cr},x}_{kT} \geq
(Nm/2)^k \geq 2^k$.

Therefore, the expectation of the total number of particles grows
at least exponentially with some exponent $\gamma > 0$. On the
other hand, from the arguments in Theorem~\ref{mtj} it follows
that for $\beta
> \beta_{\rm cr}$ the expectation of the number of particles grows exponentially
with the exponent $\lambda_0(\beta)$. Since $\lambda_0(\beta)
\downarrow 0$ as $\beta \downarrow \beta_{\rm cr}$ and $ n^{\beta
,x}_t$  depends monotonically on $\beta$, we conclude that $\gamma
= 0$. Thus we come to a contradiction with our assumption that
$S(x)$ is not identically zero.
 \qed
\\
\\
{\bf Remark 1.} It is not difficult to see that the random
variables $n^{\beta,x}_t$ can be realized on a common probability
space in such a way that $\zeta^{\beta,x}$ from Theorem~\ref{mtk}
converge almost surely, as $\beta \uparrow \beta_{\rm cr}$, to
$\eta^x$.
\\
\\
{\bf Remark 2.} Using the spectral techniques similar to those
employed above and in~\cite{CKMV1}, one can get the asymptotics of
the higher order moments of $ n^{\beta_{\rm cr},x}_t$.

\section{Appendix}
Here we prove several statements on the analytic properties of the
resolvent $R^\beta_\lambda$. We largely follow \cite{CKMV1}. The
new steps concern the asymptotics of the resolvent as $\lambda
\rightarrow 0$.
 Denote
\[
A_\lambda= v(x)R^{0}_\lambda : C_{\exp}( \mathbb{R}^d) \rightarrow
C_{\exp}( \mathbb{R}^d). \label{al1}
\]

The following lemma is similar to Lemmas 5.1 and 5.2 of
\cite{CKMV1}  (see also \cite{Va} for a similar statement for
general elliptic operators), the difference being that we now
obtain a more precise asymptotics of $ R^0_\lambda$ and
$A_\lambda$ near the origin.
\begin{lemma}
\label{ker10} Consider the operators $R^0_\lambda: C_{\exp}(
\mathbb{R}^d) \rightarrow C( \mathbb{R}^d)$ and $A_\lambda:
C_{\exp}( \mathbb{R}^d) \rightarrow C_{\exp}( \mathbb{R}^d)$.

(1) The operators $R^0_\lambda$ and $A_\lambda$  are analytic in
$\lambda \in \mathbb{C}'$.

(2) The operator $A_\lambda$ is compact for $\lambda \in
\mathbb{C}'$.

(3) For each $\varepsilon
> 0$, we have $\max(||R^0_\lambda||, ||A_\lambda||)=O(1/|\lambda|)$ as $\lambda \rightarrow
\infty$, $|{\rm arg} \lambda| \leq \pi - \varepsilon$.

(4) The operator $R^0_\lambda$ has the following asymptotic
behavior as $\lambda \rightarrow 0$, $\lambda \in \mathbb{C}' $:
\begin{equation} \label{odde}
R^0_\lambda=Q_{d}
\lambda^{\frac{d}{2}-1}+{P_{d}}(\lambda)+O(|\lambda
|^{\frac{d}{2}}),~~~d \geq3,~~~d-{\rm odd},
\end{equation}
\begin{equation} \label{evene}
R^0_\lambda=Q_{d} \lambda^{\frac{d}{2}-1} \ln ({1}/{\lambda})
+{P_{d}}(\lambda)+O(|\lambda |^{\frac{d}{2}}|\ln \lambda|),~~~d
\geq 4,~~~d-{\rm even},
\end{equation}
 where $P_{d}$ are polynomials with
coefficients which are bounded operators and
\[
Q_{d} f \equiv q_d \int_{ \mathbb{R}^d} f(x) d x,
\]
where $q_d \neq 0$.
\end{lemma}

\noindent {\bf Remark.} The term with the coefficient $Q_d$ is the
main non-analytic term of the expansion as $\lambda \rightarrow
0$.
\\
\noindent \proof Let $d$ be odd.
From  (\ref{ha2}) it follows that the kernel $R^0_\lambda (x)$ is
an analytic function of  $\lambda \in \mathbb{C}'$ and the
following estimates holds
\begin{equation} \label{estw}
|R^0_{\lambda} (x)| \leq C_d |\sqrt{\lambda}|^{d-2}
|\sqrt{\lambda} x|^{2-d},~~~|\sqrt{\lambda} x| \leq 1,
\end{equation}
\begin{equation} \label{estwT}
|R^0_{\lambda} (x)| \leq C_d |\sqrt{\lambda}|^{d-2} |e^{-\sqrt{2
\lambda}  |x|}| |\sqrt{\lambda}
x|^{\frac{1-d}{2}},~~~|\sqrt{\lambda}  x| \geq 1.
\end{equation}
Moreover, these estimate admits differentiation in $\lambda$ and
$x$ resulting in
\begin{equation} \label{estwW}
|\frac{\partial R^0_{\lambda} (x)}{\partial \lambda}| \leq C'_d
|\sqrt{\lambda}|^{d-4} |\sqrt{\lambda} x|^{2-d},~~~|\sqrt{\lambda}
x| \leq 1,
\end{equation}
\begin{equation} \label{estwWT}
|\frac{\partial R^0_{\lambda} (x)}{\partial \lambda}| \leq C'_d
|\sqrt{\lambda}|^{d-2} |e^{-\sqrt{2 \lambda}  |x|}||\sqrt{\lambda}
x|^{\frac{1-d}{2}}(\frac{1}{|\lambda|} +
\frac{|x|}{|\sqrt{\lambda}|}),~~~|\sqrt{\lambda} x| \geq 1,
\end{equation}
and
\begin{equation} \label{estDw}
|\nabla_x R^0_{\lambda} (x)| \leq C''_d |\sqrt{\lambda}|^{d-2}
|\sqrt{\lambda} x|^{2-d} \frac{1}{|x|},~~~|\sqrt{\lambda} x| \leq
1,
\end{equation}
\begin{equation} \label{estwDT}
|\nabla_x R^0_{\lambda} (x)| \leq C''_d |\sqrt{\lambda}|^{d-2}
|e^{-\sqrt{2 \lambda}  |x|}| |\sqrt{\lambda}
x|^{\frac{1-d}{2}}(\frac{1}{|x|} +
|\sqrt{\lambda}|),~~~|\sqrt{\lambda} x| \geq 1.
\end{equation}
where $C_d$, $C'_d$ and $C''_d$ are positive constants. The
estimates (\ref{estw})-(\ref{estwT}) and
(\ref{estwW})-(\ref{estwWT}) imply the analyticity  of the
operators $R^0_\lambda$ and $A_\lambda$.

The estimates (\ref{estw})-(\ref{estwT}) and
(\ref{estDw})-(\ref{estwDT}) imply that the operator $R^0_\lambda:
C_{\exp}(\mathbb{R}^d) \rightarrow C^1(\mathbb{R}^d)$ is bounded.
Then the standard Sobolev embedding theorem implies the
compactness of the operator $A_\lambda = v R^0_\lambda$ in the
space $C_{\exp}(\mathbb{R}^d)$.

In order to prove the third statement of the lemma, we observe
that the norm of $R^0_\lambda$  can be estimated by $ \int_{
\mathbb{R}^d} |{R}^{0}_\lambda (x)| d x$, which is of order
$O(1/|\lambda|)$ as
 $\lambda \rightarrow \infty$, $|{\rm arg}
\lambda| \leq \pi - \varepsilon$, due to
(\ref{estw})-(\ref{estwT}).

To prove the fourth statement, we need a more detailed asymptotics
of the kernel $R^0_\lambda (x)$ when $|\sqrt{\lambda} x|
\downarrow 0$. Namely, it follows from the properties of the
Hankel functions and (\ref{ha2}) that there are a constant $a_d
\in \mathbb{C}$ and a polynomial $b_d$ with complex coefficients
such that
\begin{equation} \label{spro}
R^0_{\lambda} (x) = |x|^{2- d}\left(a_d (\sqrt{\lambda} |x|)^{d-2}
+ b_d(\lambda |x|^2) + O(|\sqrt{\lambda} x|^{d}) \right)~~~{\rm
as}~~~|\sqrt{\lambda} x| \downarrow 0,~~~~\lambda \in \mathbb{C}'.
\end{equation}
 Combined with (\ref{estwT}) and the definition of
the space $ C_{\exp}( \mathbb{R}^d)$, this easily implies the
fourth statement of the lemma.

The case of even $d$ is similar. The main difference concerns
formulas (\ref{estw}) and  (\ref{spro}). The estimate (\ref{estw})
remains valid except the case $d = 2$, where it is replaced by
\[
|R^0_{\lambda} (x)| \leq C_2 |\ln \sqrt{\lambda}
x|,~~~|\sqrt{\lambda} x| \leq 1,
\]
while (\ref{spro}) is replaced by
\begin{equation} \label{sproM}
R^0_{\lambda} (x) = |x|^{2- d}\left(a_d (\sqrt{\lambda} |x|)^{d-2}
\ln(\sqrt{\lambda} |x|) + b_d(\lambda |x|^2) + O(|\sqrt{\lambda}
x|^{d}|\ln(\sqrt{\lambda} |x|)|) \right)
\end{equation}
 as $|\sqrt{\lambda} x| \downarrow 0$, $\lambda \in
\mathbb{C}'$. The rest of the arguments proceed as earlier, but
now employing (\ref{sproM}) instead of  (\ref{spro}). \qed
 \\

The following lemma is simply a resolvent identity.
\begin{lemma}
For $\lambda \in \mathbb{C}^{\prime}$, we have the following
relation between the meromorphic operator-valued functions
\label{lrid}
\begin{equation}  \label{aa1}
R^{\beta}_\lambda=R^0_\lambda - R^0_\lambda(I+\beta
v(x)R^0_\lambda)^{-1}(\beta v(x) R^0_\lambda)
\end{equation}
\end{lemma}

\noindent {\bf Remark.} Note that (\ref{aa1}) can be written as
\[
R^{\beta }_\lambda=R^{0}_\lambda-R^{0}_\lambda (I+\beta A_\lambda
)^{-1}(\beta v(x)R^{0}_\lambda ).  \label{aa2}
\]
From here it also follows that
\begin{equation}
R^{\beta }_\lambda =R^{0}_\lambda (I+\beta A_\lambda )^{-1},
\label{aa3}
\end{equation}
which should be understood as an identity between meromorphic in
$\lambda $ operators acting from  $C_{\exp}( \mathbb{R}^d)$ to $C(
\mathbb{R}^d)$.
\\

The kernels of the operators $I+\beta A_\lambda $, $\lambda \in
\mathbb{C}'$, are described by the following lemma.
\begin{lemma}
\label{re1} (1) The operator-valued function $(I + \beta
A_\lambda)^{-1}$ is meromorphic in $\mathbb{C}^{\prime}$. It has a
pole at $\lambda \in \mathbb{C}^{\prime}$ if and only if $\lambda$
is an eigenvalue of $L^\beta$. These poles are of the first order.

(2) Let $\lambda_i(\beta)$ be a positive eigenvalue of $L^\beta$.
There is a one-to-one correspondence between the kernel of the
operator $I + \beta A_{\lambda_i}$ and the eigenspace of the
operator $L^\beta$ corresponding to the eigenvalue $\lambda_i$.
Namely, if $(I + \beta A_{\lambda_i})h = 0$, then $\psi =-
R^0_{\lambda_i} h$ is an eigenfunction of $L^\beta$ and $h = \beta
v \psi$.

(3) If $d \geq 3$, there is a one-to-one correspondence between
the kernel of the operator $I+\beta A_0$ and solution space of the
problem
\begin{equation}
L^{\beta }(\psi )=\frac{1}{2}\Delta \psi +\beta v(x)\psi
=0,~~~\psi (x) = O(|x|^{2-d}),~~\frac{\partial{\psi}
}{\partial{r}}(x) = O(|x|^{1-d})~~as~~r =|x| \rightarrow \infty.
\label{grs}
\end{equation}
Namely, if $(I + \beta A_0)h = 0$ for $h \in  C_{\exp}(
\mathbb{R}^d)$, then $\psi =-R^0_0h$ is a solution of (\ref{grs})
and $h = \beta v \psi$.
\end{lemma}
 \proof
The operator $A_\lambda$, $\lambda \in \mathbb{C}'$, is analytic,
compact, and tends to zero as $\lambda \rightarrow +\infty$ by
Lemma~\ref{ker10}. Therefore $(I+\beta A_\lambda )^{-1}$ is
meromorphic  by the Analytic Fredholm Theorem.

 If $\lambda \in \mathbb{C}^{\prime }$ is a pole
of $(I+\beta A_\lambda)^{-1}$, then it is also a pole of the same
order of $R^{\beta }_\lambda$ as follows from (\ref{aa3}) since
the kernel of $R^{0}_\lambda $ is trivial. Therefore the pole is
simple and coincides with one of the eigenvalues $\lambda _{i}$.
Note that $\lambda $ is a pole of $(I+\beta A_\lambda)^{-1}$ if
and only if the kernel of $I+\beta A_\lambda$ is non-trivial. Let
$h\in C_{\exp}(\mathbb{R}^d)$ be such that
$||h||_{C_{\exp}(\mathbb{R}^d)}\neq 0$ and $(I+\beta
vR^{0}_\lambda )h=0$. Then $\psi :=-R^{0}_\lambda h\in
L^{2}(\mathbb{R}^{d})$ and $(\frac{1}{2}\Delta -\lambda +\beta
v)\psi =0$, that is $\psi $ is an eigenfunction of $L^{\beta }$.

Conversely, let $\psi \in L^{2}(\mathbb{R}^{d})$ be an
eigenfunction corresponding to an eigenvalue $\lambda_{i}$, that
is
\begin{equation}
(\frac{1}{2}\Delta -\lambda _{i})\psi +\beta v\psi =0. \label{bb1}
\end{equation}
Denote $h=\beta v\psi $. Then  $(\frac{1}{2}\Delta -\lambda
_{i})\psi =-h$. Thus $\psi =-R^{0}_{\lambda _{i}}h$ and
(\ref{bb1}) implies that $h$ satisfies $(I+\beta
vR^{0}_{\lambda_{i}})h=0$. Note that $h \in C^\infty(
\mathbb{R}^d)$, $h$ vanishes outside ${\rm supp} (v)$, and
therefore belongs to the kernel of $I+\beta A_{\lambda_{i}}$. This
completes the proof of the first two statements.

Similar arguments can be used to prove the last statement. If
$h\in C_{\exp}(\mathbb{R}^d)$ is such that
$||h||_{C_{\exp}(\mathbb{R}^d)}\neq 0$ and $(I+\beta A_0)h=0$,
then $h$ has compact support and  the integral operator $R^0_0$
can be applied to $h$. It is clear that $\psi :=-R^{0}_0h$
satisfies (\ref{grs}).

In order to prove that any solution of (\ref{grs}) corresponds to
an eigenvector of $I + \beta A_0$, one only needs to show that the
solution $\psi $ of the problem (\ref{grs}) can be represented in
the form $\psi =-R^{0}_0h$ with $h=\beta v\psi .$ The latter
follows from the Green formula
\[
\psi(x) =-(R^{0}_0(\beta v\psi))(x)+\int_{|y|=a}[R^{0}_0(x-y)\psi'_{r}(y)-\frac{\partial }{%
\partial r}R^{0}_0(x-y)\psi (y)]d s,\text{ \ \ \ }|x|<a,
\]
after passing to the limit as $a\rightarrow \infty$.  \qed
\\

\noindent {\bf Remark.} The relations (\ref{grs}) are an analogue
of the eigenvalue problem for zero eigenvalue and the
eigenfunction $\psi $ which does not necessarily belong to
$L^{2}(\mathbb{R}^d)$. We shall call a non-zero solution of
(\ref{grs}) a ground state.
\\


Due to the monotonicity and continuity of $\lambda =
\lambda_0(\beta)$ for $\beta > \beta_{cr}$, we can define the
inverse function
\[
\beta = \beta(\lambda): [0, \infty) \rightarrow
[\beta_{cr},\infty).
\]

We'll prove that the operator $-A_\lambda$, $\lambda
>0$, has a non-negative kernel and has a positive simple
eigenvalue such that all the other eigenvalues are smaller in
absolute value. Such an eigenvalue is called the principal
eigenvalue.

\begin{lemma} \label{ll44} The operator $-A_\lambda$, $\lambda > 0$,
has the principal eigenvalue. This eigenvalue is equal to
$1/\beta(\lambda)$ and the corresponding eigenfunction can be
taken to be positive in the interior of ${\rm supp} (v)$ and equal
to zero outside of ${\rm supp} (v)$. If $d \geq 3$, then the same
is true for the operator $-A_0$ (in particular, $\beta_{cr} > 0$).
\end{lemma}
\proof The maximum principle for the operator $(\frac{1}{2}\Delta
-\lambda )$, $\lambda
>0$, implies that the kernel of the operator $R^{0}_\lambda$,
$\lambda >0$, is negative. Thus, by (\ref{al1}), for all $y$ the
kernel of $-A_\lambda $ is positive when $x$ is in the interior of
${\rm supp} (v) $ and zero otherwise. Thus $-A_\lambda$, $\lambda
>0$, has the principal eigenvalue (see \cite{Kr}).
On the other hand, by Lemma~\ref{re1}, $1/\beta(\lambda)$ is a
positive eigenvalue of $ -A_\lambda$. Note that this is the
largest positive eigenvalue of $-A_\lambda$. Indeed, if $\mu =
1/\beta^{\prime}> 1/\beta(\lambda)$ is an eigenvalue of
$-A_\lambda$, then $\lambda$ is one of the eigenvalues $\lambda_i$
of $L^{\beta^{\prime}}$ by Lemma~\ref{re1}. Therefore, $
\lambda_i(\beta^{\prime}) = \lambda_0(\beta)$ for $\beta^{\prime}<
\beta$. This contradicts the monotonicity of $\lambda_0(\beta)$.
Hence the statement of the lemma concerning the case $\lambda
> 0$ holds.

For $d \geq 3$, the kernel of $-A_0$ is equal to $v P_d$ and has
the same properties as the kernel of $-A_\lambda$, $\lambda > 0$.
Thus $-A_0$ has the principal eigenvalue. Since $A_\lambda
\rightarrow A_0$ as $\lambda \downarrow 0$, the principal
eigenvalue $1/\beta(\lambda)$ converges to the principal
eigenvalue $\mu < \infty$ of $-A_0$. On the other hand,
$\beta(\lambda)$ is a continuous function, and therefore $\mu =
1/\beta_{cr}$, which proves the statement concerning the case
$\lambda = 0$. \qed
\\

\noindent {\bf Remark.}  Let $d \geq 3$.  Lemmas~\ref{re1} and
\ref{ll44} imply that the ground state of the operator $L^\beta$
for $\beta = \beta_{cr}$ (defined by (\ref{grs})) is defined
uniquely up to a multiplicative constant and corresponds to the
principal eigenvalue of $A_0$. If $\beta < \beta_{cr}$, then the
ground state (with $\lambda = 0$) does not exist  and the operator
$I + \beta A_0$ has a bounded inverse.
\\

%

We can finally proceed with the proof of Lemma~\ref{rbetaW}.
\\

\noindent{\it Proof of Lemma~\ref{rbetaW}.}  The analytic
properties of $R^\beta_\lambda$ follow from (\ref{aa3}) and the
corresponding properties of $(I + \beta A_\lambda)^{-1}$ which are
in turn due to Lemma \ref{re1}.

By Lemma~\ref{ker10}, the norm of $A_\lambda$ decays at infinity
when $\lambda \rightarrow \infty$, $|{\rm arg} \lambda| \leq \pi -
\varepsilon$. Therefore there is $\Lambda > 0$ such that the
operator $(I + \beta A_\lambda)^{-1}$ is bounded for
 $|{\rm arg} \lambda| \leq \pi -
\varepsilon$, $|\lambda| \geq \Lambda$. The decay of the norm of
$R^\beta_\lambda$ now follows from (\ref{aa3}) and the third part
of Lemma~\ref{ker10}.

Now let us prove (\ref{Modd}). Let $d \geq 3$ be odd. We use
(\ref{odde}) to represent $I + \beta A_\lambda$ as
\[
I + \beta A_\lambda = B + C_\lambda,
\]
where $B = I + \beta A_0$ and $C_\lambda  = \beta v(Q_{d}
\lambda^{\frac{d}{2}-1}+{P_{d}}(\lambda) - P_d(0)+O(|\lambda
|^{\frac{d}{2}}))$. Since $B$ is invertible (see the Remark
following Lemma~\ref{ll44}) and $C_\lambda \rightarrow 0$ as
$\lambda \rightarrow 0$, $\lambda \in \mathbb{C}'$, we have
\[
(I + \beta A_\lambda)^{-1} = (B + C_\lambda)^{-1} = (I - (B^{-1}
C_\lambda) + (B^{-1} C_\lambda)^2 - ...)B^{-1}
\]
for all sufficiently small $\lambda \in \mathbb{C}'$. Combining
this with (\ref{aa3}) and using (\ref{odde}), we
obtain~(\ref{Modd}).
In order to prove (\ref{Meven}), we can repeat the same arguments,
starting with (\ref{evene}) instead of~(\ref{odde}). \qed
\\
\\
{\bf \large Acknowledgements}: The author is grateful to S.
Molchanov for introducing him to this problem and valuable
discussions. While working on this article, the author was
supported by the NSF grant DMS-0854982.


\begin{thebibliography}{9}



\bibitem{ABY} Albeverio, S., Bogachev L., Yarovaya E., {\it Asymptotics
of branching symmetric random walk on the lattice with a single
source},  C. R. Acad. Sci. Paris Sér. I Math.  326  (1998),  no.
8, 975--980.

\bibitem{ABY2} Albeverio, S., Bogachev, L., Yarovaya, E., {\it Branching random walk with a single source},
Communications in difference equations (Poznan, 1998),  9--19,
Gordon and Breach, Amsterdam, 2000.



\bibitem{BY1} Bogachev L., Yarovaya E., {\it A limit theorem for a supercritical branching
 random walk on $Z^d$ with a single source},
translation in  Russian Math. Surveys  53  (1998),  no. 5,
1086--1088.

\bibitem{CaMo} Carmona, R. A.; Molchanov, S. A., {\it Stationary parabolic Anderson
model and intermittency.}  Probab. Theory Related Fields  102
(1995),  no. 4, 433--453.

\bibitem{CKMV1} Cranston M., Koralov L., Molchanov S., Vainberg B., {\it
Continuous model for homopolymers},  Journal of Functional
Analysis 256 (2009), no. 8, 2656--2696.

\bibitem{DaFl} Dawson D., Fleischmann K., {\it A super-Brownian motion with a single point catalyst},
Stochastic Process. Appl.  49  (1994),  no. 1, 3--40.

\bibitem{DFL} Dawson, D., Fleischmann K., Le Gall J., {\it Super-Brownian motions in catalytic media},
 Branching processes (Varna, 1993),  122--134, Lecture Notes in Statist., 99, Springer, New York, 1995.

\bibitem{DFLM} Dawson, D., Fleischmann, K., Li Y., Mueller C., {\it Singularity of super-Brownian
local time at a point catalyst},  Ann. Probab.  23  (1995),  no.
1, 37--55.

\bibitem{FL} Fleischmann K., Le Gall J., {\it A new approach to
the single point catalytic super-Brownian motion},  Probab. Theory
Related Fields  102  (1995),  no. 1, 63--82


\bibitem{FSh} Frechet  M., Shohat J.,
{\it A Proof of the Generalized Second-Limit Theorem in the Theory
of Probability}, Transactions of the American Mathematical Society
Vol. 33, No. 2 (Apr., 1931), pp. 533-543.

\bibitem{GK} Gartner J., Konig W., {\it The parabolic Anderson
model. Interacting stochastic systems}, 153--179, Springer,
Berlin, 2005.

\bibitem{GKM} Gartner, J., Konig, W., Molchanov, S. A., {\it Almost sure
asymptotics for the continuous parabolic Anderson model},  Probab.
Theory Related Fields  118  (2000),  no. 4, 547--573.

\bibitem{GM1} Gartner, J., Molchanov, S. A., {\it Parabolic
problems for the Anderson model. I. Intermittency and related
topics},  Comm. Math. Phys.  132  (1990),  no. 3, 613--655.

\bibitem{GM2} Gartner, J., Molchanov, S. A., {\it Parabolic problems for the Anderson
model. II. Second-order asymptotics and structure of high peaks},
Probab. Theory Related Fields  111  (1998),  no. 1, 17--55.

\bibitem{KLMS} Konig W., Lacoin H., Morters P., Sidorova. N.,
{\it A two cities theorem for the parabolic Anderson model},  Ann.
Probab. 37  (2009),  no. 1, 347--392.

\bibitem{Kr}  Krasnoselski M., {\it Positive Solutions of Operator
Equations}, Groningen, P. Noordhoff [1964].


\bibitem{Va} Vainberg B., {\it On the Short-Wave Asymptotic Behavior of
Solutions of Stationary Problems and the Asymptotic Behavior as $t
\rightarrow \infty$ of Solutions of Non-Stationary Problems},
Russian Math Surveys, 30:2, 1975, pp 1-58.

\bibitem{Ya} Yarovaya E., {Branching random walks in inhomogeneous
media} (in Russian), Moscow University Meh-Mat publication.
\end{thebibliography}
\end{document}